
\input amstex
\documentstyle{amsppt}
\pagewidth{125mm}
\pageheight{195mm}
\parindent=8mm
\frenchspacing \tenpoint

\topmatter
\title
Homological algebra of mirror symmetry
\endtitle
\author
Maxim Kontsevich
\endauthor
\affil
Max-Planck-Institut f\"ur Mathematik, Bonn\\
and University of California, Berkeley
\endaffil
\endtopmatter
\document
\nologo
\NoBlackBoxes
\TagsOnRight
\define\h{\text{H}}
\define\Q{\bold Q}
\define\C{\bold C}
\define\Z{\bold Z}

\define\R{\bold R}
\define\p{\bold P}
\define\M{\Cal M_{g,n}}
\define\Mbar{\overline{\Cal M}_{g,n}}
\define\eps{\epsilon}

Mirror Symmetry was discovered several years ago in string theory
 as a duality between families of 3-dimensional Calabi-Yau manifolds (more
precisely, complex algebraic manifolds possessing holomorphic volume elements
without zeroes).
 The name comes from the symmetry among Hodge numbers. For dual Calabi-Yau
manifolds $V,\,\,W$ of dimension $n$ (not necessarily equal to $3$) one has
$$\text{dim }\,\h^p(V, \Omega^q)=\text{dim }\,\h^{n-p}(W, \Omega^q)\,\,\,.$$

Physicists conjectured that conformal field theories associated with mirror
varieties are equivalent. Mathematically, MS is considered now as a relation
between numbers of rational curves on  such a manifold and Taylor coefficients
of periods of Hodge structures considered as functions on the moduli space of
complex structures on a mirror manifold. Recently it has been realized that one
can make predictions for numbers of curves of positive genera and also on
Calabi-Yau manifolds of arbitrary dimensions.

We will not describe here the complicated history of the subject and will not
mention many beautiful contsructions, examples and conjectures motivated by MS.
On the contrary, we want to give an outlook of the story in general terms and
propose a conceptual framework for a possible explanation of
 the
 mirror phenomenon. We will restrict ourselves to a half of MS considering
it as a relation between symplectic structures on one side and complex
structures on another side. Actually, we will deal only with a half of this
half, ignoring the holomorphic anomaly effects (see [BCOV]) in the symplectic
part (A-model) and the polarization of Hodge structures in the complex part
(B-model).
 For an introduction to Mirror Symmetry we recommend [M] and [Y].

At the moment there are only few completely solid statements, essentially
because there were no universal  definition of the ``number of curves'' for a
long time.

\subheading{ Comparision of symplectic and complex geometry}

We start with a recollection of well-known facts concerning  symplectic and
complex manifolds.

Let $V$ be a compact smooth $2n$-dimensional manifold.

1.S. A symplectic structure on $V$ is given by a reduction of the structure
group $GL(2n,\R)$ of the tangent bundle $T_V$ to the subgroup $Sp\,(2n,\R)$
satisfying certain integrability conditions (the associated $2$-form $\omega$
is closed, or, equivalently the associated Poisson bracket  on smooth functions
satisfies the Jacobi identity).

1.C. A complex structure on $V$ is given by a reduction  of the structure group
$GL(2n,\R)$ of the tangent bundle $T_V$ to the subgroup $GL(n,\C)$ satisfying
certain integrability conditions (the Newlander-Nirenberg theorem).

Notice that both groups $Sp\,(2n,\R)$ and $GL(n,\C)$ are homotopy equivalent to
$U(n)$. Thus, they have the same algebra of characteristic classes generated by
Chern classes $c_i\in \h^{2i}(BU(n),\Z),\,\,\,1\le i \le n$.

Basic examples of compact symplectic or complex manifolds are complex
projective algebraic  manifolds endowed with the pull-back of the Fubini-Studi
K\"ahler form on the projective space.

2.S. First-order deformations of symplectic structures on $V$ are in one-to-one
correspondence with $\h^2(V,\R)$. The deformation theory is unobstructed and
the  local moduli space of symplectic structures on $V$ can be identified with
a domain in the affine space $\h^2(V,\R)$ via map $\omega\mapsto[\omega]\in
\h^2(V,\R)$ (J.~Moser).

2.C. First-order deformations of complex structures on $V$ near a fixed one are
in one-to-one correspondence with $\h^1(V,T_V^{hol})$ where $T_V^{hol}$ denotes
the sheaf of holomorphic vector fields on $V$ (Kodaira-Spencer theory). If
$c_1(V)=0$ and $V$ admits a K\"ahler structure then the deformation theory of
$V$ is unobstructed and the local moduli space can be identified with a domain
in the affine space $\h^1(V,T_V^{hol})$ (the Bogomolov-Tian-Todorov theorem).

The following two facts concern only complex manifolds.

3.C. For a complex manifold $V$ admitting a K\"ahler structure there is a pure
Hodge structure on the singular cohomology groups: $$\h^k(V,\Z)\otimes\C\simeq
\bigoplus_{p+q=k}\h^q(V,\Omega^p)\,\,\,.$$

4.C. With a complex algebraic manifold $V$ one can associate the abelian
category $Coh(V)$ of coherent sheaves on $V$ and the triangulated category
$\Cal D^b(Coh(V))$ (the bounded derived  category).

Our aim in this talk is to propose candidates for 4.S and 5.S in the context of
symplectic geometry. The mirror symmetry should be a correspondence (partially
defined and multiple valued) between symplectic and complex manifolds (both
with $c_1=0$) identifying structures 2-4.

To get a feeling of what is going on it is instructive to look at a simplest
case of the mirror symmetry, which is already highly non-trivial.

\subheading{Two-dimensional tori (after R.~Dijkgraaf)}

Let $\Sigma$ be a complex elliptic curve and $p_1,\dots,p_{2g-2}$ are pairwise
distinct points of $\Sigma$, where $g\ge 2$ is an integer.
  We consider holomorphic maps $\phi: C@>>>\Sigma$ from compact connected
smooth  complex curves $C$ to $\Sigma$, which have  only one double
ramification point over each point $p_i\in \Sigma$ and no other ramification
points. By the Hurwitz formula the genus of $C$ is equal to $g$. The set
$X_{g}(d)$ of equivalence classes of such maps of degree $d\ge 1$ is finite,
and for each $\phi: C@>>>\Sigma$ its automorphism group
$$\text{Aut}\,\,(\phi):=\{ f:C@>>>C\,\,|\,\,\,\phi\circ f=\phi\}$$
is finite.  For $g\ge 2$ we introduce the generating series in one variable $q$
as follows:
$$F_g(q):=\dsize\sum_{d\ge 1}\left(\sum_{[\phi]\in X_{d}(d)} \frac
{1}{|\,\text{Aut}\,\,(\phi)\,\,|}\,\right){q^d}\,\in\Q[[q]]\,\,\,.$$

The following statement is now rigorously established due to efforts of several
people (R.~Dijkgraaf, M.~Douglas, D.~Zagier, M.~Kaneko):
$$F_g\in\Q[E_2,E_4,E_6]\,\,\,,$$
where $E_k$ are the classical Eisenstein series,
 $$E_k=1-\frac{2k}{B_k}\sum_{n\ge 1} \left(\sum_{a|n}a^{k-1}\right)
q^n\,\,\,.$$
 $E_k$ is a modular form of weight $k$ for even $k\ge $ and $E_2$ is \it not
\rm a modular form.\newline
 Here $B_2=1/6,\,\,B_4=-1/30,\,\,B_6=1/42,\,\,\dots$ are Bernoulli numbers. If
one associates with $E_k,\,\,k=2,4,6$, the degree $k$, then $F_g$ has degree
$6g-6$.

One can regard  $\Sigma$ as a symplectic 2-dimensional manifold
 $(S^1\times S^1, \omega)$
 with the area $\int_{\Sigma}\omega$ equal to $-\text{log}\,(q),\,\,\,0<q<1$,
and interpret weights  $q^d$ of ramified coverings as
$$\text{exp}\,(-\text{area of }C\text{ with respect to the pullback of
}\omega)\,\,\,.$$
Mirror symmetry in this example is the claim that the generating function for
certain  invariants of symplectic structures on $S^1\times S^1$ is a ``nice''
function on the moduli space of complex structures on $S^1\times S^1$. The
two-dimensional torus is a self-dual manifold for MS.

Notice that the standard local coordinate $q=\,\,\text{exp}\,(2\pi i
\tau)\,,\text{ Im}\,\tau>0$ on the moduli space of elliptic curves,
 $$\text{Elliptic curve}\,=\C/\left(\Z+\Z\tau\right)\,\,,$$
 can be written as
$$q=\,\,\text{exp}\,\left(2\pi i\frac{\int_{\gamma_1}\Omega}
{\int_{\gamma_0}\Omega}\right)\,\,$$
where $\gamma_0,\,\gamma_1$ are two generators of $\h_1(\text{elliptic
curve},\Z)$ and $\Omega$ is a non-zero holomorphic $1$-form.

\subheading{Quintic $3$-folds (after [COGP])}

Here we describe the first famous prediction of physicists.  Let $V$ be a
non-singular hypersurface in complex projective space $\p^4$ given by an
equation \newline
$Q(x_1,\dots,x_5)=0$ of degree $5$ in $5$ homogeneous variables
$(x_1:\dots:x_5)$. This complex manifold carries a  top degree holomorphic
differential form of non-degenerate at all points (a holomorphic volume
element):
 $$\Omega=\frac{1}{dQ} \sum_{i=1}^{5}(-1)^i x_i dx_1\wedge\dots
\wedge\widehat{dx_i}\wedge\dots \wedge dx_5\,\,.$$

H.~Clemens conjectured that smooth rational curves on a generic quintic
$3$-fold are isolated. Recently it was checked up to degree 9.
It is natural to count the number $N_d$ of rational curves on $V$ of fixed
degree $d$. In fact, there are singular rational curves on $V$ of degree $5$,
and one has to take them into account as well. At the end of  the next section
we will propose an algebro-geometric formula for the ``physical'' number of
curves on $V$, both smooth and singular,  without assuming the validity of the
Clemens conjecture.

The mirror symmetry prediction is the following.
 First of all, we define the virtual number of curves of degree $d$ as
$$N_d^{virt}:=\sum_{k|d}\frac{1}{k^3}N_{d/k}\,\in \,\Q\,\,\,.$$
The reason for this fomula is that in string theory one counts not just curves
in $V$ but maps from rational curves to $V$. Any map $\p^1@>>>V$ of positive
degree is the composition of a rational map $\p^1@>>>\p^1$ and of an embedding
$\p^1\hookrightarrow V$. It was argued first by Aspinwall and Morrison [AM]
that the  factor associated with multiple coverings of degree $k$ should be
equal to $1/k^3$.

 The complete generating function for rational curves is
$$F(t):=\frac{5}{6}t^3+\sum_{d\ge 1} N_d^{virt} \text{exp}\,(dt)\,\,.$$
The first summand here represents the contribution of maps of degree $0$ (i.e.,
maps to a point of $V$).

 On the mirror side we consider functions
$$\psi_0(z)=\sum_{n=0}^{\infty} \frac{(5n)!}{(n!)^5} z^n$$
$$\psi_1(z)=\,\,\text{log}\,\,z\cdot \psi_0(z)+5\sum_{n=1}^{\infty}
\frac{(5n)!}{(n!)^5} \left(\sum_{k=n+1}^{5n} \frac{1}{k}\right)z^n$$
$$\psi_2(z)=\frac{1}{2}\,(\,\text{log}\,z)^2\cdot \psi_0(z)+\dots$$
$$\psi_3(z)=\frac{1}{6}\,(\,\text{log}\,z)^3\cdot \psi_0(z)+\dots$$
which are solutions of the linear differential equation
$$\left(\left(z\frac{d}{dz}\right)^4-5z(5z\frac{d}{dz}+1)(5z\frac{d}{dz}+2)(5z\frac{d}{dz}+3)(5z\frac{d}{dz}+4)\right)\psi(z)=0\,\,\,.$$
More precisely,
$$\sum_{i=0}^3
\psi_i(z)\eps^i+O(\eps^4)=\sum_{n=0}^{\infty}\frac{(1+5\eps)\,(2+5\eps)\,\dots\,(5n+5\eps)}{\left((1+\eps)\,(2+\eps)\,\dots\,(n+\eps)\right)^5} \,\,z^{n+\eps}\,\,\,.$$

Functions $\psi_i(z)$ are periods $\int_{\gamma_i}\omega$ of the Calabi-Yau
manifold $W=W(z)$ which is a resolution of singularities of the following
singular variety:
$$\{(x_1:x_2:x_3:x_4:x_5)|\,\,x_1^5+x_2^5+x_3^5+x_4^5+x_5^5=z\,x_1 x_2 x_3 x_4
x_5\}/(\Z/5\Z)^3\,\,.$$
Here the group $(\Z/5\Z)^3$ is the group of diagonal matrices preserving the
equation,
$$\{\,\,\text{diag}(\xi_1,\dots,\xi_5)\,| \,\,\,\xi_i^5=1,\,\,\prod_{i=1}^5
\xi_i=1\,\}/\{\,\xi\,\text{Id}\,\,|\,\,\xi^5=1\,\}$$
 and $\gamma_i$ are certain singular homology classes with complex
coefficients.  Family of varieties $W(z)$ depending on $1$ parameter is mirror
dual to a universal family of smooth quintic $3$-folds depending on  $101$
parameters.

The predictions of physicists is that
$$F\left(\frac{\psi_1}{\psi_0}\right)=\frac{5}{2}\,\cdot\,\frac{\psi_1\psi_2-\psi_0\psi_3}{\psi_0^2}\,\,\,.$$

One of miracles in this formula is that
$$\text{exp}\,\left(\frac{\psi_1}{\psi_0}\right)\in \Z[[z]]\,\,.$$

Also, numbers $N_d$ computed via the mirror prediction are positive integers.

It is interesting that the contribution of \it individual \rm non-parametrized
rational curves on $3$-dimensional Calabi-Yau manifolds is connected with
variations of mixed Hodge structures in a fashion analogous to the mirror
symmetry predictions. Namely, according to the Aspinwall and Morrison formula
[AM] we have the following generating function:
$$F(t)=\sum_{d\ge 1} \frac{\text{exp}\,(dt)}{d^3}\,\,\,.$$
We introduce functions $\psi_*$:
$$\psi_0(z)=1\,,\,\,\,\psi_1(z)=\text{log}\,z\,,\,\,\,\psi_2(z)=\frac{1}{2}\,(\,\text{log}\,z)^2\,,\,\,\,\psi_3(z)=\text{Li}_3(z)\,:=\sum_{d\ge 1}\frac{z^d}{d^3}$$
which are solutions of the linear differential equation
$$\frac{d}{dz}\left(\frac{1-z}{z}\left(z\frac{d}{dz}\right)^3
\,\psi(z)\right)=0\,\,.$$
Functions $F$ and $\psi_*$ are related by the evident formula
  $$F\left(\frac{\psi_1}{\psi_0}\right)=\frac{\psi_3}{\psi_0}\,\,\,.$$

\subheading{Gromov-Witten invariants}

We describe here a not yet completely constructed theory which has potentially
wider domain of applications than mirror symmetry.
  It is based on pioneering ideas of M.~Gromov [G] on the role of
$\overline{\partial}$-equations in symplectic geometry, and  certain  physical
intuition proposed by E.~Witten [W1], [W2]. There are many   evidences that the
following picture from [KM] is correct.

Let $(V,\omega)$ be a closed symplectic manifold, $\beta\in \h_2(V,\Z)$ be a
homology class and $g,n\ge 0$
 be integers satisfying the inequality $2-2g-n<0$.
 Gromov-Witten classes
$$I_{g,n;\beta}\in \h_D\left(\Mbar(\C)\times V^n;\Q\right)$$
are homology classes with rational coefficients of degree
  $$D=D(g,n,\beta)=(\text{dim}\,V-6)(1-g)+2n+2\int_{\beta} c_1(T_V)\,\,\,.$$
Here $\Mbar$ denotes the Deligne-Mumford compactification of the moduli stack
of smooth connected algebraic curves of genus $g$ with $n$ marked points.
Recall that an algebraic curve $C$ with marked points $p_1,\dots,p_n$ is called
\it stable \rm if
\roster
\item all singular points of $C$ are ordinary double points,
\item marked points $p_i$ are pairwise distinct and smooth, $p_i\in C^
{smooth}$,
\item the group of automorphisms of $(C:p_1,\dots,p_n)$ is finite, or,
equivalently, the Euler characteristic of each connected component of
$C^{smooth}\setminus\{p_1,\dots,p_n\}$ is negative.
\endroster
 The arithmetic genus  of  stable curve $C$ is defined by the formula
$$g_a(C):=\text{ dim}\,\h^1(C,\Cal O)=\chi(C^{smooth})/2-1\,\,\,.$$
 The stack $\Mbar$ is the moduli  stack of stable marked curves of arithmetic
genus $g$ with $n$ marked points.
The associated coarse moduli space $\Mbar(\C)$  is a compact complex orbifold.

   One expects that $I_{g,n;\beta}$ is  invariant under  continuous
deformations of the symplectic structure on $V$.

 Analogously, we expect that the Gromov-Witten invariants can be defined for
non-singular projective algebraic varieties over arbitrary fields and they take
values in the Chow groups with rational coefficients instead of the singular
homology groups.

Intuitively, the geometrical meaning of  Gromov-Witten classes in the
symplectic case can be described as follows. Let us choose an almost-complex
structure on $V$ compatible in the evident way with the symplectic form
$\omega$. Notice that the space of almost-complex structures compatible with
the fixed $\omega$ is contractible.
 Denote by $X_{g,n}(V,\beta)$ the space of equivalence classes of
$(C;x_1,\dots,x_n;\phi)$ where
 $C$ is a smooth complex curve of genus $g$ with pairwise distinct marked
points $x_i$, and $\phi:C@>>>V$ is a pseudo-holomorphic map (i.e., a solution
of the Cauchy-Riemann equation $\overline{\partial}\phi=0$) such that the image
of the fundamental class of $C$ equal to $\beta$.  There is a natural map from
$X_{g,n}(V,\beta)$ to $\M(\C)\times V^n$
 associating with $(C;x_*;\phi)$ the equivalence class of $(C;x_*)$ and the
sequence of points $(\phi(x_1),\dots,\phi(x_n))$. One can show easily that the
dimension of the space $X_{g,n}(V,\beta)$ at each point is bigger than or equal
to $D(g,n,\beta)$. Also, under  appropriate genericity assumptions,
$\,\,X_{g,n}(V,\beta)$ is \it smooth \rm of dimension $D(g,n,\beta)$.
 We want to define
$I_{g,n;\beta}$ as the image of the fundamental class of a  compactification
$\overline{X_{g,n}(V,\beta)}$. The problem here is to find a correct
compactification and to define the ``fundamental class'' if there are
components of dimensions bigger than  expected.  Also one has to prove that
classes $I_{g,n;\beta}$ do not depend on the choice of an almost-complex
structure.

There are now two approaches to the rigorous construction of Gromov-Witten
classes. First one is due to Y.~Ruan and G.~Tian [RT] and it suffices   for the
genus zero case. This construction works only for so called semi-positive
manifolds (including Fano and Calabi-Yau manifolds), but it gives classes with
integral coefficients.  The idea of this construction is to perturb generically
$\overline{\partial}$-equations and check that there are no strata of
dimension larger than  $D(0,n,\beta)$ in  Gromov's compactification of the
space of pseudo-holomorphic curves.
 In fact, Ruan and Tian define not GW-classes but the number  of maps from a
fixed complex curve to $V$ satisfying general incidence conditions (counted
with signs). Using algebraic results on the structure of $\h^*(\overline{\Cal
M}_{0,n})$ it is possible to reconstruct whole genus-zero part of Gromov-Witten
classes (see [KM]).
Another construction [K1] is based on a new compactification of the moduli
space of maps and should work, presumably, for all genera, for all symplectic
manifolds and also for all non-singular projective varieties over arbitrary
fields. At least, one can produce now purely algebro-geometric definitions of
genus-zero Gromov-Witten invariants in the case of complete intersections in
 projective spaces. Its advantage is that it will not use any general position
argument, and its weak point is the lack of  control on  integrality of arising
classes.

As an example we give a definition of ``numbers of rational curves'' on a
quintic $3$-fold. Denote by $\overline{\Cal M}_{0,0}(\p^4,d)$ the moduli stack
of equivalence classes of maps
 $\phi:C@>>>\p^4$ where $C$ is a connected curve of arithmetic genus zero and
having only odrinary double points as  singular points (i.e., $C$ is a tree of
rational curves) such that each irreducible component of $C$ mapping to a point
has at least $3$ singular points. The parameter $d,\,\,d\ge 1\,,\,\,$ denotes
the degree of the image of the fundamental class $[C]$ in
 $\h_2(\p^4,\Z)\simeq \Z$.
  It is proven in [K1] that $\overline{\Cal M}_{0,0}(\p^4,d)$ is a smooth
proper algebraic stack of finite type. The set of its complex points is a
compact complex  orbifold of dimension $5d+1$.

We define a vector bundle $\Cal E_d$ of rank $5d+1$ over $\overline{\Cal
M}_{0,0}(\p^4,d)$. The fiber of $\Cal E_d$ at $\phi:C@>>>\p^4$
 is equal to $\h^0(C,\phi^*\Cal O (5))$.
  Notice that if a quintic $3$-fold $V$ is given by an equation $Q$ of degree
$5$ in $5$ variables, $Q\in \Gamma(\p^4,\Cal O (5))$,
  then there is an associated section $\widetilde{Q}$ of $\Cal E_d$ whose
zeroes are exactly maps  into $V$. In general, there are connected components
of the set of zeroes of $\widetilde{Q}$
 of positive dimensions arising from  multiple covering maps to rational curves
in $V$. Nevertheless, we define the ``virtual'' number of curves by the formula
$$N_d^{virt}:=\int\limits_{\overline{\Cal M}_{0,0}(\p^4,d)} c_{5d+1}(\Cal
E_d)\,\,\,.$$

Here the integral is understood in the orbifold sense. Thus, the numbers
 $N_d^{virt}\in \Q$  in general are not integers.
We are  sure that our formula will give the same numbers as physicists predict.
This formula was checked up to degree $4$. Also we obtained in [K1] a closed
expression for the generating function for the numbers $N_d^{virt}$. Hence the
mirror prediction in the quintic case is reduced to an explicit identity.

There is an extension of the definition  above to the case of complete
intersections in toric varieties and for  counting of higher genus curves in
flag varieties.

\subheading{Axioms}

The system of axioms formulated in [KM] is a formalization of what physicists
call 2-dimensional topological field theory coupled with gravity (see [W2]). We
reproduce here only one of axioms from [KM] which is the basic one. Other
axioms encode more evident properties of Gromov-Witten classes, like the
invariance under permutation of indices, etc.

It will be convenient to associate with the class
$I_{g,n;\beta}$ a linear map
$$J_{g,n;\beta}:(\h^*(V,\Q))^{\otimes n}@>>>\h^*(\Mbar,\Q)$$
using the K\"unneth formula and the Poincar\'e duality.
  A \it splitting \rm axiom describes the restriction of Gromov-Witten classes
to boundary divisors of $\Mbar$. Namely, for $g_1,g_2\ge 0$ and $n_1,n_2\ge 0$
such that
$$g_1+g_2=g,\,\,\,n_1+n_2=n+2,\,\,\,2-2g_i-n_i<0 \text{ for }i=1,2$$
  there exists a natural inclusion $i_{g_*,n_*}:\overline{\Cal
M}_{g_1,n_1}\times \overline{\Cal M}_{g_2,n_2}\hookrightarrow\Mbar $
$$(C_1;x_1,\dots,x_{n_1})\times (C_2;y_1,\dots,y_{n_2})\longmapsto
(C_1\bigcup_{x_1=y_1}C_2;x_2,\dots,x_{n_1},y_2,\dots,y_{n_2})\,\,.$$

The following diagram should be commutative:

$$\CD
\h^*(V)^{\otimes n}@>\simeq>>
\h^*(V)^{\otimes (n_1-1)}\otimes \h^*(V)^{\otimes (n_2-1)} \\
@V J_{g,n;\beta} VV @VV\otimes\Delta V \\
\h^*(\Mbar) @.
\h^*(V)^{\otimes n}\otimes\h^*(V)^{\otimes n_2}\\
@V(i_{g_*,n_*})^* VV @VV\sum\limits_{\beta_1+\beta_2=\beta}
 J_{g_1,n_1;\beta_1}\otimes J_{g_2,n_2;\beta_2}V \\
\h^*(\overline{\Cal M}_{g_1,n_1}\times\overline{\Cal
M}_{g_2,n_2})@>\simeq>\text{K\"unneth}>
\h^*(\overline{\Cal M}_{g_1,n_1})\otimes \h^*(\overline{\Cal M}_{g_2,n_2})
\endCD
$$

Here all cohomology are taken with coefficients in $\Q$ and $\Delta$ denotes
the Poincar\'e dual to the fundamental class
 of the diagonal $V\subset V\times V$.
 The geometric meaning of this axiom is clear:
 a map $\phi$ of the glued curve $C$ from the image of $i_{g_*,n_*}$ is the
same as two maps $\phi_1,\,\phi_2$ from $C_1$ and $C_2$ with
$\phi_1(x_1)=\phi_2(y_1)$.

The splitting axiom in the case $g_1=g_2=0$ was checked by Y.~Ruan and G.~Tian
for semi-positive manifolds and by me for complete intersections using the
stable map approach.

\subheading{Associativity equation}

For a compact symplectic manifold $(V,\omega)$ denote by $\Cal H:=\oplus_k
\h^k(V,\C)$ the total cohomology space of $V$ considered as a $\Z$-graded
vector space (super vector space) and also as a complex supermanifold. It means
that the underlying topological space of $\Cal H$ is $\h^{even}(V,\C)$ and
functions $\Cal O(\Cal H)$ on $\Cal H$ are holomorphic functions on
$\h^{even}(V,\C)$ with values in the exterior algebra generated by
$\left(\h^{odd}(V,\C)\right)^*$.

Using Gromov-Witten classes for genus zero we define the following function
(pre-potential) on $\Cal H$:
$$\Phi(\gamma):=\sum\limits_{\beta\in \h_2(V,\Z)} e^{-\int_{\beta}
\omega}\sum\limits_{n\ge 3}
\frac{1}{n!}\int\limits_{I_{0,n;\beta}}1_{\overline{\Cal
M}_{0,n}}\otimes\gamma\otimes \dots\otimes\gamma\,.$$
Here $\gamma$ denotes an even element of $\Cal H \otimes \Lambda$ where
$\Lambda$ is an arbitrary auxilary supercommutative algebra (as usual in the
theory of supermanifolds). The element  $1_{\overline{\Cal M}_{0,n}}$ is the
identity in the cohomology ring of  $\overline{\Cal M}_{0,n}$.

\proclaim{Conjecture} The series $\Phi$ is absolutely convergent in a
neighbourhood  $\Cal U$ of $0$ in $\Cal H$, if the cohomology class $[\omega
]\in \h^2(V,\R)$ is sufficiently positive.
\endproclaim

Without assuming the validity of this conjecture we can work not over the field
$\C$ but over the field of fractions of the semigroup  ring $\Q[B]$ where $B$
is the semigroup generated by classes $\beta$ such that $\int_{\beta}
\omega'\ge 0$ for  all symplectic froms $\omega'$ close to $\omega$. Other
homology classes are excluded because they cannot be represented by
pseudo-holomorphic curves.

The function $\Phi$ in its definition domain $\Cal U$ satisfies a system of
non-linear differential equations of the third order (due to R.~Dijkgraaf,
E.~Verlinde, H.~Verlinde, and E.~Witten, see [W2]). Let us choose a basis $x_i$
of the space $\Cal H$ and denote by $x^i$
the corresponding  coordinate system on $\Cal H$. Denote by $(g_{ij})$ the
matrix of the Poincar\'e pairing, $\,g_{ij}:=\int_V  x_i\wedge x_j $ ,  and by
$(g^{ij})$  the inverse matrix.
For all $i,j,k,l$, we have (modulo appropriate sign corrections  for odd-degree
classes):
$$
\dsize\sum_{m,m'}\frac{\partial^3 \Phi}{\partial x^i
\partial x^j\partial x^m} g^{mm'}
\frac{\partial^3 \Phi}{\partial x^{k}
\partial x^l\partial x^{m'}}=
\dsize\sum_{m,m'}\frac{\partial^3 \Phi}{\partial x^i
\partial x^k\partial x^m} g^{mm'}
\frac{\partial^3 \Phi}{\partial x^{j}
\partial x^{l}\partial x^{m'}}
$$
This equation can be reformulated as the condition of associativity of the
algebra given by the structure constants $A_{ij}^{k}:=\sum_{k'}
g^{kk'}\partial_{ijk'}\Phi$. In invariant terms it means that $\Phi$ defines a
supercommutative associative multiplication on the tangent bundle to $\Cal H$
(the quantum cohomology ring).

The associativity equation follows from the splitting axiom and from a certain
linear relation among components of the compactification divisor of
$\overline{\Cal M}_{0,n}$.
Denote by $D_S$ for $S\subset \{1,\dots,n\},\,\,2\le \# S\le n-2,\,\,$
 the divisor in $\overline{\Cal M}_{0,n}$ which is the closure of the moduli of
stable curves $(C;p_1,\dots,p_n)$ consisting of two irreducible components
$C_1,C_2$ such that $p_i\in C_1$ for $i\in S$ and $p_i\in C_2$ for $i\notin S$.
\proclaim{Lemma} We have the following identity in $\h^2(\overline{\Cal
M}_{0,n},\Z)$
$$\dsize\sum\Sb S:1,2\in S
\\\,\,\,\,3,4\notin S\endSb [D_S]=\sum\Sb S:1,3\in S
\\\,\,\,\,2,4\notin S\endSb [D_S]\,\,.$$
\endproclaim

Both sides in the equality above are pullbacks under the forgetful map\newline
 $\overline{\Cal M}_{0,n}@>>>\overline{\Cal M}_{0,4}$ of points
$D_{\{1,2\}}$, $D_{\{1,3\}}\in \overline{\Cal M}_{0,4}\simeq \p^1$.
 It is clear that any two points on $\p^1$ are rationally equivalent as
divisors.
\qed

Conversely, one can show using the splitting axiom that one can reconstruct the
whole system of genus zero GW-classes starting from $\Phi$. The equation of the
associativity is a necessary and sufficient condition for the existence of such
 reconstruction (see [KM]).

The associativity equation was studied by B.Dubrovin [D]. He discovered that it
is a completely integrable system in many cases (but not for CY manifolds). For
example, for $V\simeq \p^2$ the associativity equation is equivalent to the
Painlev\'e VI equation.
 It means that via a simple recursion formula we can compute the number of
rational curves of degree $d$ in the projective plane passing through generic
$3d-1$ points.

Notice that by dimensional reasons, the associativity equation is an empty
condition for $3$-dimensional Calabi-Yau manifolds, because the virtual
dimension of the space of rational curves is zero, curves do not intersect each
other and the degeneration argument is unapplicable.

Let us introduce a connection on the tangent bundle $T_{\Cal U}$ by the formula
$\nabla={\nabla_0}_{|\Cal U}+A$, where $\nabla_0$ is the standard connection of
the affine space $\Cal H$.
 The associativity equation implies the flatness of $\nabla$.

\subheading{Variations of Hodge structures}

Suppose that $c_1(V)=0$, and $V$ carries at least one integrable complex
structure compatible with $\omega$ such that $\h^{2,0}(V)=0$. For any such
complex structure we have a Hodge decomposition $\oplus\h^k(V,\C)=\oplus
\h^{p,q}$. We expect that all cycles $I_{g,n;\beta}$ are Hodge cycles  of
(complex) dimension equal to $(n+\text{dim}_{\C}V-3)$.  It follows that the
restriction of $\nabla$ to the convergence domain of the series $\Phi$ in the
second cohomology group:
 $$\Cal U^{cl}:=\Cal U\cap \h^2(V,\C)\subset\h^2(V,\C)=\h^{1,1}$$  maps
$\h^{p,q}\otimes \Cal O(\Cal U^{cl})$ to $\h^{p+1,q+1}\otimes\Omega^1(\Cal
U^{cl})$. We call $\Cal U^{cl}$ the classical moduli space because it is
locally isomorphic to a complexification of the moduli space of symplectic
structures on $V$.

We introduce filtrations $\oplus_{p\le p_0}\h^{p,q}$ on trivial bundles over
$\Cal U^{cl}$ with fibers equal to
 $\oplus_{p-q\text { is fixed}} \h^{p,q}$. Hence we have flat connections and
filtrations on holomorphic vector bundles over a complex manifold satisfying
the Griffiths transversality conditions.
 We call such data a \it complex \rm variation of pure Hodge structures.
One can prove  by using formal arguments with Hodge-Tate groups that the
equivalence classes of such complex variations of pure Hodge structures do not
change under deformations of the complex structure on $V$.

For general symplectic manifolds $V$ with $c_1(V)=0$ we can consider just the
two trivial vector bundles $\h^{ev}$ and $\h^{odd}$ on $\Cal U^{cl}:=\Cal
U\cap\h^2(V,\C)\subset \Cal H$ endowed with the flat connection induced from
$\nabla$  and the filtration  by subbundles $\oplus_{k\le k_0} \h^k$.

 Algebro-geometric complex variations of pure Hodge structures are defined as
  subquotients of variations of pure Hodge structures on cohomology groups of
complex projective algebraic manifolds depending algebraically on parameters.

\proclaim{Mirror Conjecture} Complex variations of pure Hodge structures
constructed using Gromov-Witten invariants of symplectic manifold $V$ as above
are locally equivalent to  algebro-geometric variations. \endproclaim

 In almost all known examples such variations of Hodge structures should be
locally equivalent to variations of Hodge structures on total cohomology
bundles  of a mirror family of complex  manifolds with $c_1=0$. Exceptions come
mostly from CY-manifolds $V$ such that $\,\,\text{dim}\,\h^1(V,T_V)=0$, i.e.
rigid manifolds.
 In this case dual manifolds with rotated Hodge diamond could not exist,
because $\,\,\text{dim}\,\h^1(W,\Omega^1_W)\ne 0$ for K\"ahler manifolds.
Physicists proposed as candidates certain substructures of Hodge structures on
cohomology groups of Fano varieties (=algebraic manifolds with an ample
anti-canonical bundle).
 Also, calculations of numbers of curves on projective spaces suggest
that in general there exists some relation between the pre-potential of \it non
\rm Calabi-Yau manifolds and algebro-geometric
variations of Hodge structures.

In the case of a quintic $V$ in $\p^4$ the function $\Phi$ is the sum of two
terms: the contribution of  maps to points of $V$ and the contribution of
rational curves in $V$ (and their multiple covers). We introduce coordinates
$t^i,\,\,i=0,1,2,3$ in one-dimensional spaces $\h^{i,i}(V)$ and odd coordinates
$\xi^j,\eta^j,\,\,\,j=1,\dots,102\,\,$ in $\h^3(V,\C)$. In these coordinates we
have (modulo adding a polynomial of degree 2)
$$\dsize\Phi(t^i,\xi^j,\eta^j)=\frac{5}{6}\sum_{i+j+k=3}t^i t^j t^k +
t^0\sum_j \xi^j\eta^j+\sum_{d\ge 1}N_d^{virt}\,\text{exp}\,( dt^1)\,\,.$$
One can deduce an  example from [COGP] from this formula.

 The flat coordinates $x^i$ on the moduli space of complex structures on dual
manifolds are equal to the ratios of periods
$\left(\int_{\gamma_i}\Omega\right)/\left(\int_{\gamma_0}\Omega\right)$ where
$\Omega$ is a holomorphic volume element on the mirror manifold $W$ and
$\gamma_i$ are elements of $H_*(W,\C)$ locally constant with respect to the
Gauss-Manin connection on the homology bundle.

There exists a generalization of the mirror correspondence to higher genera.
 First of all, the dimension formula for degrees of Gromov-Witten classes shows
that one can expect a non-negative dimension for the space of genus $g$ curves
for Calabi-Yau varieties $V$ only in the following cases:
\roster
\item $g=0$ and an arbitrary dimension $n:=\text{dim}\,V$ (this is what we have
described right now),
\item $g=1$ and arbitrary $n$,
\item $g\ge 2$ and $n\le 3$.
\endroster

The Harvard group of physicists in  the remarkable paper [BCOV] proposed a
procedure (``quantum Kodaira-Spencer theory'') giving numbers of curves for
cases $g=1$ or $n=3$. It relates GW-invariants with
 certain structures on the moduli of dual varieties which are more complicated
than just  variations of Hodge structures and are not understood mathematically
yet. The example of R.~Dijkgraaf (elliptic curves) is a $1$-dimensional version
of this theory.\newline
\newline

In the rest of this talk we give an outline of a program relating Mirror
Symmetry  to general structures of Homological Algebra. The central ingredient
here is a fundamental construction of K.~Fukaya based on ideas of S.~Donaldson,
A.~Floer and G.~Segal.
\subheading{Extended moduli spaces}

When we restrict the flat bundle $T_{\Cal U}$ to the subspace $\Cal
U^{cl}=\h^2(V,\C)$, much of information will be lost. It seems very reasonable
to extend the moduli space of symplectic structures to the whole domain $\Cal
U$ in $\Cal H$ in which the potential $\Phi$ is defined. Hence the
  tangent space to the extended moduli space at classical points $\Cal U^{cl}$
should be equal to $\Cal H=\oplus \h^k$.

  Now we want to construct an extended moduli space $\Cal{M}$ for a complex
Calabi-Yau $W$ containing the ordinary moduli space  of complex structures on
$W$. The natural candidate for the tangent bundle to $\Cal{M}$ at classical
points  $\Cal M ^{cl}:=\text{Moduli}\,(W)$ should be equal to the direct sum
$\oplus \h^p(W,\wedge^q T_W)$. The problem of constructing $\Cal M$ was already
dicussed by E.~Witten (see  [W3]).

  We anticipate that $\oplus \h^p(W,\bigwedge^q T_W)$ can be interpreted as the
total Hochschild
 cohomology of the sheaf $\Cal{O}_W$ of algebras of holomorphic functions on
$W$.

For an algebra $A/k$ over a field its Hochschild cohomology
$\h\h^*(A)=\h^*(A,A)$  is defined as $\text{Ext}^*_{A-mod-A}(A,A)$.
 The second Hochschild cohomology   $\h\h^2(A)$  classifies
 infinitesimal deformations of $A$.
 Notice that each $A$-bimodule $M$ defines a functor from the category of
$A$-modules into itself:
 $$M\otimes\Sb A\endSb:\,\,\,A-mod\,@>>>\,A-mod\,,\,\,\,\,N\mapsto M\otimes\Sb
A \endSb N$$
and $A$ corresponds to the identity functor $\text{Id}_{A-mod}$.

Analogously, we define the Hochschild cohomology of the structure  sheaf $\Cal
O_W$ of a scheme $W$ over $k$ (or of an analytic space)  as the global
$\text{Ext}$-functor
 $$\h\h^*(\Cal O_W):=\,\text{Ext}^*_{W\times W}\left(\delta_*(\Cal
O_W),\delta_*(\Cal O_W)\right)$$
where $\delta:W\hookrightarrow W\times W$ is the diagonal embedding.
  Another definition of the Hochschild cohomology for algebraic varieties (in
fact, equivalent to the ours) was proposed by M.~Gerstenhaber and S.~D.~Schack
[GS].
  The following fact proven in hidden form in [GS] seems to be new in algebraic
geometry:
\proclaim{Theorem} For smooth (and not necessarily compact) variety $W$ over a
field of characteristic zero there is a canonical isomorphism
  $$\h\h^n(\Cal O_W)\simeq\bigoplus_{k+l=n} \h^k(W,{\bigwedge}^l T_W)\,\,.$$
\endproclaim

For smooth  $W$ the second Hochschild cohomology $\h\h^2(\Cal O_W)$ splits
into the direct sum of ordinary first-order deformations $\h^1(W,T_W)$,
non-commutative deformations $\h^0(W,\bigwedge^2 T_W)$ of the sheaf $\Cal O_W$
of associative algebras (global Poisson brackets on $W$), and a little bit more
misterious piece $\h^2(W,\Cal O_W)$. This third part can be interpreted as
locally trivial first-order deformations of the sheaf of abelian categories of
$\Cal O_W$-modules.

In the next section we will propose an interpretation of the total Hochschild
cohomology as the tangent space to  ``extended moduli space''
 $\Cal M$ containing the classical moduli space $\Cal M^{cl}$ as a part.

\subheading{$A_{\infty}$-algebras and categories
}

$A_{\infty}$-algebras were introduced by J.Stasheff in 1963 (see [S]). Let
$A=\oplus A^k$ be a $\Z$-graded vector space. The structure of
$A_{\infty}$-algebra on $A$ is an infinite sequence of linear maps
$m_k:A^{\otimes k}@>>>A,\,\,\,k\ge 1,\,\,\,\text{deg}\,m_k=2-k$  satisfying the
(higher) associativity conditions:
\roster
\item $m_1^2=0$, (we can consider $m_1$ as a differential and $(A,m_1)$ as a
complex),
\item $m_1(m_2(a\otimes b)=m_2(m_1(a)\otimes b)\pm m_2(a\otimes m_1(b))$,
($m_2$ is a morphism of complexes),
\item $m_3(m_1(a)\otimes b\otimes c)\pm m_3(a\otimes m_1(b)\otimes c)\pm
m_3(a\otimes b\otimes m_1(c))\pm m_1(m_3(a\otimes b \otimes c)=$
$$=m_2(m_2(a\otimes b)\otimes c)-m_2(a\otimes m_2(b\otimes c)),\text{ ($m_2$ is
associative up to homotopy)},$$
\item  and so on...
\endroster

In one sentence one can define the $A_{\infty}$-algebra structure on $A$ as a
co-derivation in the graded sense $d,\,\,\,d^2=0$ of degree $1$ on the co-free
co-associative algebra without a co-unit
co-generated by the $\Z$-graded vector space $A[1],\,\,\,\,\,A[1]^k:=A^{k+1}$.

A morphism of $A_{\infty}$-algebras (from $A$ to $B$) is an infinite collection
of linear maps $A ^{\otimes k}@>>>B\,,\,\,\,k\ge 1$ satisfying some equations
analogous to the defining equations for  individual $A_{\infty}$-algebras.
 In terms of co-free co-algebras it is the same as a differential graded
homomorphism. A homotopy equivalence of $A_{\infty}$-algebras is a morphism
whose linear part induces an  isomorphism of cohomology groups with respect to
the differential $m_1$.

In general, $A_{\infty}$-algebras are closely related to  differential graded
algebras. Namely, a dg-algebra is the same as an $A_{\infty}$-algebra with
$m_3=m_4=\dots=0$.  Conversely, for an $A_{\infty}$-algebra $A$ one can
construct using the bar-construction a differential graded algebra $B$
 homotopy equivalent to $A$.

An additive category over a field $k$ is a category $C$ with finite direct sums
such that  all sets of morphisms $\text{Hom}_C(X,Y)$ are  endowed with
structure of vector spaces over $k$ and where the composition of morphisms is a
bilinear map. In a sense, one can approximate additive categories by algebras
of endomorphisms of their objects.
 Analogously, one can define  differential graded category as an additive
category with the structure of complexes on  $\text{Hom}_C(X,Y)$ such that the
composition is a morphism of complexes.

An
$A_{\infty}$-category $C$ is a collection of objects and $\Z$-graded spaces of
morphisms $\text{Hom}_C(X,Y)$ for each two objects endowed with higher
compositions of morphisms satisfying relations parallel to the defining
relations of $A_{\infty}$-algebras.  We require the existence of identity
morphisms $Id_X\in\text{Hom}_C(X,X)$ which obey the usual properties of
identity for composition $m_2$ and vanish under substitution in other (higher)
compositions. We can also require the existence of finite direct sums in $C$ in
an obvious sense.
 Notice that $C$ is not a category in general, because the composition of
morphisms is not associative. Nevertheless one can construct an additive
category $H(C)$ from $C$ with the same class of objects by defining new
$\Z$-graded spaces of morphisms as
$$\text{Hom}_{H(C)}(X,Y):=\frac{\text{Ker}\,\,(m_1:\text{Hom}_C^0(X,Y)@>>>\text{Hom}_C^1(X,Y)\,)}{\text{Im}\,\,(m_1:\text{Hom}_C^{-1}(X,Y)@>>>\text{Hom}_C^0(X,Y)\,)}\,\,\,\,.$$

There exists a generalization of Hochschild cohomology to the case of
$A_{\infty}$-algebras. The meaning of $\h\h^*(A)$ for $*>0$ is the space of
equivalence classes of first-order deformations of $A_{\infty}$-structure on
$A$ over $\Z$-graded bases.   We hope that there exists an appropriate version
of Hochschild cohomology for some good class of $A_{\infty}$-categories as
well.

One can show under some mild assumptions that if the formal $\Z$-graded
 moduli space $\Cal M$ of $A_{\infty}$-categories is smooth than there exists
the canonical
 structure of an associative commutative algebra on the tangent bundle to $\Cal
M$.
 In the case of a an $A_{\infty}$-category consisting just of one object $X$
with morphisms $\text{Hom}_C(X,X)$ consisting of an associative algebra $A$ in
degreee $0$
  the product in Hochschild cohomology (i.e. in the tangent space to $\Cal M$)
$$\h\h^*(A):=\,\text{Ext}^*_{A-mod-A}(A,A)$$
 coincides with the usual Yoneda composition of $Ext$-groups.

\subheading{Triangulated categories}

One of fundamental tools in homological algebra is the triangulated category
$\Cal D(C)$ associated with an abelian category $C$ satisfying certain
conditions (J.-L.~Verdier, see [V]). A triangulated category is an additive
category endowed with a shift functor  and a class of so-called exact
triangles, obeying a complicated list of axioms. For example, for $C$ equal to
the category of $A$-modules, where $A$ is an associative algebra, the category
$\Cal D(C)$ is equivalent to the category whose objects are complexes of free
$A$-modules and morphisms are equal to homotopy classes of differential graded
morphisms of degree $0$:
$$\text{Hom}_{D(C)}(X,Y):=\h^0(\bigoplus_k\prod_j \text{Hom}_C(X^j,
Y^{j+k}))\,\,\,.$$

The bounded derived category $\Cal D^b(C)$ is the  full subcategory of $\Cal
D(C)$ consisting of complexes of $A$-modules with non-vanishing cohomology
groups only in finitely many degrees.

The shift functors at the level of objects just shifts degree of complexes:
$X@>>>X[n],\,\,\,X[1]^k=X^{k+n}$ and $(X[n])[m]=X[n+m],\,\,\,X[0]=X$.

 We will not describe Verdier's axiomatics of exact triangles here because it
does not look completely satisfactory, although it was generally adopted and
widely used. Certain improvement of axioms was proposed by A.~Bondal and
M.~Kapranov in [BK]. The main ingredient in their definition is the notion of a
twisted complex in a differential graded category.

We can extend  the construction of [BK]  to the case of
$A_{\infty}$-categories.  We assume that an $A_{\infty}$-category $C$ is
endowed with  shift functors such that
$$ \text{Hom}_C(X[i], Y[j])=\text{Hom}_C(X, Y)\,[j-i]\,\,\,.$$

 By definition, a (one-sided) \it twisted complex \rm  is a family
$(X^{(i)})_{i \in \Z}$ of objects of an $A_{\infty}$-category $C$ such that
$X^{(i)}=0$ for almost all $i$ together with a collection of morphisms
$d_{ij}\in \text{Hom}_C(X^{(i)},
X^{(j)})^{j-i}$ for $i<j$ obeying a generalization of the Maurer-Cartan
equation:
 $$\text{for fixed }\,i,\,j \,\,\,\,\,\,\,\,\,\,\dsize\sum\Sb
k;i_0,\dots,i_k\\i_0=i,i_k=j\endSb
m_k(f_{i_0,i_1},\dots,f_{i_{k-1},i_k})=0\,\,\,. $$
We define the $\Z$-graded space of morphisms between twisted complexes $X$ and
$Y$ as $$\bigoplus_{k,j} \text{Hom}_C(X^{(j)}, Y^{(j+k)})[-k]\,\,\,.$$
Using higher compositions in $C$ one can define the structure of an
$A_{\infty}$-category on $\{\text{twisted complexes of }C\}$. Any  higher
composition of morphisms of twisted complexes is defined as the sum over all
possible products which one can imagine.

One can check without difficulties that the derived category
$$\Cal D^b(C):=H(\text{twisted complexes of }C)$$
 satisfies Verdier axioms for triangulated category.

\subheading{Fukaya's $A_{\infty}$-category}

 In this section we describe a remarkable contsruction of Kenji Fukaya [F] with
few minor modifications.

Let $V$ be a closed symplectic manifold with $c_1(T_V)=0$.

Denote by $L V$ the  space of pairs $(x,L)$ where $x$ is a point of $V$ and
$L$ is a Lagrangian subspace in $T_x V$.  The space $ L V$ is fibered over $V$
with fibers equal to Lagrangian Grassmanians. Thus the fundamental group of the
fibers is isomorphic to $\Z$.

The condition  on $V$  posed above guarantees that there exists
 a $\Z$-covering $\widetilde{ L V}$ of $  L V$ inducing a universal cover of
each fiber. Let us fix $\widetilde{ L V}$.

Objects of Fukaya's category $F(V)$ are Lagrangian submanifolds $\Cal L\subset
V$ endowed with a continuous lift of the evident map $\Cal L@>>>  L V$ to a map
$\Cal L@>>>\widetilde{ L V}$. In fact, it is only a first approximation to
right objects, see remarks after the defintion.
  For subvarieties $\Cal L_1, \Cal L_2$ intersecting each other transversally
at a point $x\in V$ and endowed with lifts to $\widetilde{ L V}$, we can define
the Maslov index $\mu_x(\Cal L_1,\Cal L_2)\in \Z$.
Notice that
$$\mu_x(\Cal L_1,\Cal L_2)+\mu_x(\Cal L_2,\Cal
L_1)=n:=\frac{1}{2}\text{dim}(V)\,\,.$$

K.~Fukaya defines the space of morphisms $Mor_F(\Cal L_1,\Cal L_2)$ only if
$\Cal L_1, \Cal L_2$ intersect transversally:
$$\text{Hom}_F(\Cal L_1,\Cal L_2):=\C^{\Cal L_1\cap\Cal  L_2}$$
with $\Z$-grading coming from the Maslov index.

The differential in $\text{Hom}_F(\Cal L_1,\Cal L_2)$ is a version of Floer's
differential. Its matrix coefficient associated with two intersection points
$p_1,p_2\in \Cal L_1\cap \Cal L_2$ is defined as
$$\sum_{\phi:D^2@>>>V} \pm \text{exp}(-\text{area of }D^2)$$
where $\phi:D^2@>>>V$ is a pseudo-holomorphic map from the standard disc
$D^2:=\{z\,|\,\,\,|z|\le 1\}\in \C$ to $V$ such that
 $\phi(-1)=p_1,\,\,\phi(+1)=p_2$ and
$$\phi(\text{exp}(i \alpha))\in \Cal L_1\text{ for } 0<\alpha<\pi,
\,\,\,\phi(\text{exp}(i \alpha))\in \Cal L_2\text{ for }
\pi<\alpha<2\pi\,\,\,.$$
More precisely, we consider \it equivalence classes \rm of maps $\phi$ modulo
the action of the group of holomorphic automorphisms of $D^2$ stabilizing
points $1$ and $-1$:
$$\R^{\times}_+\subset \,PSL(2,\R)=\text{ Aut }(D^2)\,\,\,.$$
The area of $D^2$ with respect to the pullback of $\omega$ depends only on the
homotopy type of $\phi\in \pi_2(V,\Cal L_1\cup \Cal L_2)$. One expects that for
sufficiently large $\omega$ the infinite series is absolutely convergent.

The sign $\,\,\pm\,\,$ in the definition of Floer differential comes from a
natural orientation of the space of pseudo-holomorphic maps. One expects that
there will be finitely many such maps for a generic almost-complex structure on
$V$ if $\mu_{p_2}-\mu_{p_1}=1$. Presumably, one can develop a general technique
of stable maps for surfaces with boundaries or extend Ruan-Tian's methods.

Analogously, one can
 define higher order compositions using zero-dimensional components of spaces
of equivalence classes modulo $PSL(2,\R)$-action of maps $\phi$ from the
standard disc $D^2$  to $V$ with the boundary $\phi(\partial D^2)$ sitting in a
union of Lagrangian subvarieties.
 More precisely, if $\Cal L_1,\dots,\Cal  L_{k+1}$ are Lagrangian submanifolds
intersecting each other transversally and $p_j\in \Cal L_j\cap \Cal
L_{j+1},\,\,\,j=1,\dots,k$ are chosen intersection points, then we define the
composition of corresponding base elements in spaces of morphisms as
$$m_k(p_1,\dots,p_k):=\sum\Sb\phi:D^2@>>>V,\,q\in \Cal L_1\cap\Cal L_{k+1}\\
0=\alpha_0<\alpha_1<\dots<\alpha_k<\alpha_{k+1}=2\pi\endSb \pm\,\,
\text{exp}\,\left(-\int_{D^2}\phi^*\omega\right)\,q \in \text{Hom}_F(\Cal
L_1,\Cal L_{k+1})$$
where $\phi(\text {exp} (i\alpha))\in \Cal L_j$ for
$\alpha_{j-1}<\alpha<\alpha_j$
 and $\phi(\text {exp} (i\alpha_j))=p_j$ for $j=1,\dots,k+1;\,\,p_{k+1}:=q$.
 Again, we expect that there exist only finitely many equivalence classes
modulo the action of $\text{Aff}(\R)=\text{Stab}_{1\in D^2}\subset PSL(2,\R)$
of such maps in each homotopy class if $\mu_q=\mu_{p_1}+\dots+\mu_{p_k}+2-k$
and the infinite series in the definition of $m_k$ converges absolutely.

K.~Fukaya claims that the identities of $A_{\infty}$-category follow from
considerations analogous to the proof of the associativity equations in the
case of rational curves. Also he claims that it is possible to extract an
actual $A_{\infty}$-category with compositions of all morphisms using an
appropriate notion of a ``generic'' Lagrangian manifold. In particular, it is
possible to restore
the identity morphisms. The main idea is that two Lagrangian submanifolds
obtained one from another by a Hamiltonian flow are equivalent with respect to
the Floer cohomology.

There is an extension of Fukaya's category. We can consider  pairs consisting
of a Lagrangian submanifold $\Cal L$ and a unitary local system $\Cal E$ on
$\Cal L$ as objects of new $A_{\infty}$-category. Morphism spaces will be
defined as
$$\text{Hom}_F((\Cal L_1, \Cal E_1),(\Cal L_2,\Cal E_2)):=\bigoplus_{p\in \Cal
L_1\cap \Cal L_2} \text{Hom}\,({\Cal E_1}_{|p},{\Cal E_2}_{|p})\,\,\,.$$

In the definiton of higher composition we add a new factor to each term equal
to the trace of the composition of holonomy maps along the boundary of $D^2$.
Unitarity in this definition is an obligatory condition, otherwise the series
defining higher compositions will be inevitably divergent.

It seems that there are further possible extensions of Fukaya's
$A_{\infty}$-category. One can consider as objects  Lagrangian foliations,
families of Lagrangian submanifolds parametrized by closed oriented manifolds
etc.

\subheading{Homological Mirror Conjecture}

We propose here a conjecture in slightly vague form which should imply the
``numerical'' Mirror conjecture. Let $(V,\omega)$ be a $2n$-dimensional
symplectic manifold with $c_1(V)=0$ and $W$ be a dual $n$-dimensional complex
algebraic manifold.

\it The derived category constructed from the Fukaya category $F(V)$ (or a
suitably enlarged one) is equivalent to the derived category of coherent
sheaves on a complex algebraic variety $W$. \rm

More precisely, we expect that there is an embedding of $\Cal D^b(F(V))$ as a
full triangulated subcategory into $\Cal D^b (Coh(W))$.
We have following evidences for that.
\roster
\item By the general philosophy, $\,\,A_{\infty}$-deformations of first order
of $F(V)$ should
 correspond to $Ext$-groups in a  category of functors $F(V)@>>>F(V)$.
 The natural candidate for such a category is $F(V\times V)$ where the
symplectic structure on $V\times V$ is $(\omega,-\omega)$. The diagonal
$V_{diag}\subset V\times V$ is a Lagrangian submanifold and it corresponds to
the identity functor. By a version of Floer's theorem (see [F]) there is a
canonical isomorphism between the Floer cohomology $\h^*(\text{Hom}_{F(V\times
V)}(V_{diag},V_{diag}))$ and the ordinary topological cohomology $\h^*(V,\C)$.
 The Yoneda product on the Floer cohomology considered as $Ext$-groups arises
from holomorphic maps from $D^2$ with $3$ marked points on $\partial D^2$ to
$V\times V$ with a boundary on $V_{diag}$.
 Such maps are the same as holomorphic maps to $V$ from the 2-dimensional
sphere
 $S^2\simeq \C P^2$ with $3$ marked points. Thus, it seems very reasonable to
expect that we will get exactly the quantum cohomology product on $\h^*(V)$.
  We expect that the equivalence of derived categories will imply numerical
predictions.
\item Lagrangian varieties (and local systems on it) form a natural class of
local boundary conditions for the  A-model in  topological open string theory.
Also, holomorphic vector bundles form local boundary conditions for the
B-model (E.~Witten [W4]). Physicists beleive that the whole string theories  on
dual varieties are equivalent. Thus, we want to say that topological open
string theory is more or less the same as a triangulated category.
\item Both categories $\Cal D^b(F(V))$ and $\Cal D^b (Coh(W))$ possess a
duality:
 a functorial isomorphism $(\text{Hom}(X,Y))^*\simeq \text{Hom}(Y,X[n])$.
 On the algebro-geometric side it is Serre duality. For Fukaya's category
 the definition of compositions is cyclically symmetric. The duality follows
from this symmetry and from the idenitity $\mu_x(\Cal L_1,\Cal L_2)+\mu_x(\Cal
L_2,\Cal L_1)=n$.
 We developed some time ago a theory of $A_{\infty}$-algebras with duality in
[K2] and proposed a combinatorial construction of cohomology classes of the
moduli spaces of smooth curves $\Cal M_{g,n}$ based on such algebras. This
construction has a generalization to an $A_{\infty}$-category with a duality.
Thus, we expect that the Gromov-Witten
invariants could be defined in a general purely algebraic situation. We still
do not know what is  missed in algebraic structures and how to define classes
with values in $H^*(\Mbar,\C)$.
\endroster

A mirror complex manifold $W$ is  usually not unique. For example, one cannot
distinguish the $B$-models on dual abelian varieties $A,\,A'$. It is compatible
with our picture because the derived categories of coherent sheaves  on $A$ and
$A'$ are equivalent via the Fourier-Mukai transform.
 Also, the $B$-models on birationally equivalent Calabi-Yau manifolds $W,\,W'$
are beleived to be isomorphic. In all known examples the Hodge structures on
total cohomology depend only on a birational type. Thus, we expect that the
derived categories of coherent sheaves on $W$ and on $W'$ are equivalent.

Our conjecture, if it is true, will unveil the mystery of Mirror Symmetry.
 The numerical predictions mean that two elements of an uncountable set
(formal power series with integral coefficients) coincide. Our homological
conjecture is equivalent to the coincidence in a \it countable \rm set
(connected components of the ``moduli space of triangulated categories'',
whatever it means).

In the last section we show what  our program looks like in the simplest case
of Mirror Symmetry.

\subheading{Two-dimensional tori: a return}

Let $\Sigma$ be the standard  flat two-dimensional torus $S^1\times S^1$
endowed with a symplectic form $\omega$ proportional to the standard volume
element.
 Let $\Cal L_1,\Cal L_2,\Cal L_3$ be three simple closed geodesics from
pairwise different homology classes and
$$p_1\in \Cal L_1\cap\Cal  L_2,\,\,\,p_2\in \Cal L_2\cap\Cal
L_3,\,\,\,p_3=q\in \Cal L_1\cap \Cal L_3$$
 be three intersection points. We will compute now the tensor coefficient of
the composition $m_2$ corresponding to the base vectors $p_1,p_2,p_3$.
 Each map $\phi$ from $D^2$ to $\Sigma$ can be lifted to a map
$\tilde{\phi}$ from $D^2$ to $\R^2=\text{the universal covering space of }
\Sigma$.
  The preimages of  circles $\Cal L_i$ on  $\R^2$ form three families of
parallel straight lines.  Thus the images of lifted maps $\tilde{\phi}$ are
triangles  with sides on these lines. It is easy to see that the equivalence
classes of triangles modulo the action of $\Z^2=\pi_1(\Sigma)$ are labeled by
terms of an arithmetic progression (the lengths of sides of triangles sitting
on the pullback of $\Cal L_1$). The areas of triangles are proportional to the
squares of elements of this progression.
 The tensor element of composition $m_2$ can be written naturally as
$$\sum_{n\in \Z} \text{ exp }(-(an+b)^2)$$
for some real parameters $a\ne 0$ and $b$,
which is a value of the classical $\theta-$function. The associativity equation
is equivalent to the standard bilinear identity for $\theta$-functions.
 It is well-known that $\theta$-functions form natural bases of spaces of
global sections of line bundles over elliptic curves.

It seems very plausible that  the triangulated category constructed from  the
Fukaya category $F(\Sigma)$ enlarged using unitary local systems is equivalent
to  the bounded derived category of coherent sheaves on the elliptic curve with
the real modular parameter $\tau:=\text{exp}\,(-\text{ area of }(\Sigma))$.

\enddocument